\documentclass[11pt]{article}
\input{diagxy}
\usepackage{amsmath} 
\usepackage{amsthm}
\usepackage{verbatim}
\usepackage{enumerate}
\usepackage{xspace}
\theoremstyle{definition}
\newtheorem{ihom}{Lemma}

\newtheorem{lift}[ihom]{Lemma}
\newtheorem{closed}[ihom]{Theorem}
\newtheorem{closed2cat}[ihom]{Corollary}
\newcommand{\pb}{\ensuremath{\textnormal{Cat}_\textnormal{pb}}\xspace}
\newcommand{\bpb}{\ensuremath{\textbf{Cat}_\textbf{pb}}\xspace}
\newcommand{\ih}[1]{\ensuremath{\textnormal{[A,#1]}_\textnormal{pb}}\xspace}
\newcommand{\un}[1]{\ensuremath{[\textnormal{A},\textnormal{#1} \times \textnormal{A}]_\textnormal{pb}}\xspace}
\newcommand{\vcat}[1]{\ensuremath{\textnormal{#1-} \textnormal{Cat} \xspace}}
\newcommand{\uno}{\ensuremath{[\textnormal{A},\textnormal{B} \times \textnormal{A}]}\xspace}
\newcommand{\iho}{\ensuremath{\textnormal{[A,B]}}\xspace}
\title{The category of categories with pullbacks is cartesian closed}
\author{John Bourke \\
email:{\tt johnb@maths.usyd.edu.au}}
\date{10/04/2009}
\begin{document}
\def\xypic{\hbox{\rm\Xy-pic}}
\label{firstpage}
\maketitle
\begin{abstract}
We show that the category of categories with pullbacks and pullback preserving functors is cartesian closed.
\end{abstract}
Consider the category, \pb, whose objects are categories with pullbacks and whose morphisms are pullback preserving functors.  Our aim is to show this category is cartesian closed.\\
 Firstly observe that \pb has products.  Given A,B of \pb, the cartesian product in Cat, $\textnormal{A} \times \textnormal{B}$, has pullbacks, constructed pointwise.  As pullbacks in $\textnormal{A} \times \textnormal{B}$ are pointwise, the projections from $\textnormal{A} \times \textnormal{B}$ preserve them.  It is straightforward to check the universal property, and so the cartesian product in \pb is just the ordinary cartesian product of categories.\\
Upon observing that \pb has products, we proceed to show that it is cartesian closed.  To do so is to provide a right adjoint  to the functor $\xy
{\ar^{\textnormal{-} \times \textnormal{A}}(0,0)*+{\pb};(600,0)*+{\pb}};
\endxy$ for each object $\textnormal{A}$ of \pb.  In keeping with convention the right adjoint is denoted by \ih{-} and referred to as the internal hom.
Given B of \pb, we define the internal hom \ih{B} to be the category whose objects are pullback preserving functors from A to B  and whose morphisms are cartesian natural transformations between such functors.
(Recall that a natural transformation is cartesian if its naturality squares are pullback squares).\\
The first thing we should verify is that \ih{$\textnormal{B}$} is actually an object of \pb, which is to say that it has pullbacks.
\begin{ihom}
\ih{$\textnormal{B}$} has pullbacks.
\end{ihom}
\begin{proof}
Given pullback preserving functors F,G and H from A to B, and a pair of cartesian natural transformations
$\xy
{\ar@{=>}^{t}(0,0)*+{F};(300,0)*+{H}};
\endxy$ and 
$\xy
{\ar@{=>}^{u}(0,0)*+{G};(300,0)*+{H}};
\endxy$ we must construct the pullback in \ih{B}.  As the category B has pullbacks, pullbacks in the ordinary functor category \iho exist and are constructed pointwise.  Consider the pullback in \iho, 
$$\xy 
(0,0)*+{P}="p";(300,0)*+{F}="f"; (0,-300)*+{G}="g"; (300,-300)*+{H}="h";
{\ar@{=>}^{r} "p";"f"}; {\ar@{=>}_{s} "p";"g"}; {\ar@{=>}^{t} "f";"h"}; {\ar@{=>}^{u}"g";"h"}; 
(50,-100)*{};(100,-100)*{} **\dir{-};
(100,-50)*{};(100,-100)*{} **\dir{-};
(500,-300)*{.};
\endxy$$  We will show that this is the pullback in \ih{B}.  To do so we must firstly show that this square lives in \ih{B}, which is to say that P preserves pullbacks and that $r$ and $s$ are cartesian.
\begin{itemize}
\item P preserves pullbacks:\\
(This is just a special case of the general fact that the limit of a diagram of $\Phi$-continuous functors is always $\Phi$-continuous, for a type of limit $\Phi$.  We give an explicit proof below).\\
Given a pullback square in A, say
$$\xy 
(0,0)*+{a}="p";(300,0)*+{c}="f"; (0,-300)*+{b}="g"; (300,-300)*+{d}="h";
{\ar_{\alpha} "p";"g"}; {\ar^{\beta} "p";"f"}; {\ar^{\theta} "g";"h"}; {\ar^{\phi}"f";"h"}; 
(50,-100)*{};(100,-100)*{} **\dir{-};
(100,-50)*{};(100,-100)*{} **\dir{-};
(500,-300)*{,};
\endxy$$ we must show that the square
$$\xy 
(0,0)*+{Pa}="p";(0,-300)*+{Pb}="f"; (400,00)*+{Pc}="g"; (400,-300)*+{Pd}="h";
{\ar_{P\alpha} "p";"f"}; {\ar^{P\beta} "p";"g"}; {\ar^{P\theta} "f";"h"}; {\ar^{P\phi}"g";"h"}; 
\endxy$$
is a pullback too.
In the composite square
$$\xy 
(0,0)*+{Pa}="p";(400,0)*+{Fa}="f"; (0,-300)*+{Pb}="g"; (400,-300)*+{Fb}="h";
{\ar^{r_{a}} "p";"f"}; {\ar_{P\alpha} "p";"g"}; {\ar_{F\alpha} "f";"h"}; {\ar_{r_{b}}"g";"h"}; 
(50,-120)*{};(100,-120)*{} **\dir{-};
(100,-70)*{};(100,-120)*{} **\dir{-};
(800,0)*+{Fc}="Fc";(800,-300)*+{Fd}="Fd";
{\ar^{F\beta} "f";"Fc"}; {\ar_{F\theta} "h";"Fd"}; {\ar^{F\phi} "Fc";"Fd"}; 
(450,-120)*{};(500,-120)*{} **\dir{-};
(500,-70)*{};(500,-120)*{} **\dir{-};
\endxy$$ both smaller squares are pullbacks, the left square as $r$ is cartesian, the right hand square as F preserves pullbacks.  Thus the composite square is a pullback.  The equations $F\beta \circ r_{a} = r_{c} \circ P\beta$ and $F\theta \circ r_{b} = r_{d} \circ P\theta$ hold by naturality of $r$ and so we may rewrite the above square as 
$$\xy 
(0,0)*+{Pa}="p";(400,0)*+{Pc}="f"; (0,-300)*+{Pb}="g"; (400,-300)*+{Pd}="h";
{\ar^{P\beta} "p";"f"}; {\ar_{P\alpha} "p";"g"}; {\ar_{P\phi} "f";"h"}; {\ar_{P\theta}"g";"h"}; 
(800,0)*+{Fc}="Fc";(800,-300)*+{Fd}="Fd";
{\ar^{r_{c}} "f";"Fc"}; {\ar_{r_{d}} "h";"Fd"}; {\ar^{F\phi} "Fc";"Fd"}; 
(450,-120)*{};(500,-120)*{} **\dir{-};
(500,-70)*{};(500,-120)*{} **\dir{-};
(950,-300)*{.};
\endxy$$
The right hand square in the composite is a pullback as $r$ is cartesian.  Therefore as the composite square is a pullback, the left hand square must be a pullback too.
\item$r$,$s$ are cartesian:
We shall consider the case of r.\\
Given a morphism $\xy
{\ar^{\alpha}(0,0)*+{a};(300,0)*+{b}};
\endxy$ of A we must show that the square 
$$\xy 
(0,0)*+{Pa}="p";(0,-300)*+{Pb}="f"; (400,00)*+{Fa}="g"; (400,-300)*+{Fb}="h";
{\ar_{P\alpha} "p";"f"}; {\ar^{r_{a}} "p";"g"}; {\ar^{F\alpha} "g";"h"}; {\ar^{r_{b}}"f";"h"}; 
\endxy$$ is a pullback.  Both squares in the composite
$$\xy 
(0,0)*+{Pa}="p";(400,0)*+{Ga}="f"; (0,-300)*+{Fa}="g"; (400,-300)*+{Ha}="h";
{\ar^{s_{a}} "p";"f"}; {\ar_{r_{a}} "p";"g"}; {\ar_{u_{a}} "f";"h"}; {\ar_{t_{a}}"g";"h"}; 
(50,-120)*{};(100,-120)*{} **\dir{-};
(100,-70)*{};(100,-120)*{} **\dir{-};
(800,0)*+{Gb}="Fc";(800,-300)*+{Hb}="Fd";
{\ar^{G\alpha} "f";"Fc"}; {\ar_{H\alpha} "h";"Fd"}; {\ar^{u_{b}} "Fc";"Fd"}; 
(450,-120)*{};(500,-120)*{} **\dir{-};
(500,-70)*{};(500,-120)*{} **\dir{-};
\endxy$$
are pullbacks: the left hand square because pullbacks in the functor category are constructed pointwise; the right hand square because $u$ is cartesian.  Therefore the composite square is a pullback.  By naturality of $t$ and $s$ we may rewrite this composite as 
$$\xy 
(0,0)*+{Pa}="p";(400,0)*+{Pb}="f"; (0,-300)*+{Fa}="g"; (400,-300)*+{Fb}="h";
{\ar^{P\alpha} "p";"f"}; {\ar_{r_{a}} "p";"g"}; {\ar_{r_{b}} "f";"h"}; {\ar_{F\alpha}"g";"h"}; 
(800,0)*+{Gb}="Fc";(800,-300)*+{Hb}="Fd";
{\ar^{s_{b}} "f";"Fc"}; {\ar_{t_{b}} "h";"Fd"}; {\ar^{u_{b}} "Fc";"Fd"}; 
(450,-120)*{};(500,-120)*{} **\dir{-};
(500,-70)*{};(500,-120)*{} **\dir{-};
(950,-300)*{.};
\endxy$$
The right hand square is a pullback, and so as the composite is, it follows that the left hand square is a pullback.
\end{itemize}
Thus the square $$\xy 
(0,0)*+{P}="p";(300,0)*+{F}="f"; (0,-300)*+{G}="g"; (300,-300)*+{H}="h";
{\ar@{=>}^{r} "p";"f"}; {\ar@{=>}_{s} "p";"g"}; {\ar@{=>}^{t} "f";"h"}; {\ar@{=>}^{u}"g";"h"}; 
(50,-100)*{};(100,-100)*{} **\dir{-};
(100,-50)*{};(100,-100)*{} **\dir{-};
(400,-300)*{.};
\endxy$$ lies in \ih{B}.  \\
Its universal property is easily checked upon noting that if $r_{1}$ and $r_{2}$ are vertically composable natural transformations such that $r_{2}$ is cartesian and $r_{2} \circ r_{1}$ is cartesian, then $r_{1}$ is cartesian.
\end{proof} 
Given a morphism $\xy
{\ar^{F}(0,5)*+{B};(300,5)*+{C}};
\endxy$ of \pb, there is an induced pullback preserving functor
$$\xy
{\ar^{\ih{F}}(0,0)*+{\ih{B}};(600,0)*+{\ih{C}}};
\endxy$$ given by composition with F, and so we obtain an endofunctor 
$$\xy
{\ar^{\ih{-}}(0,0)*+{\pb};(600,0)*+{\pb}};
(800,0)*{.};
\endxy$$ We shall show that \ih{-} is right adjoint to $\textnormal{-} \times \textnormal{A}$ by providing a unit and counit satisfying the triangle equations.  These may be lifted to \pb directly from the case of Cat.
In the case of the cartesian closedness of Cat, the unit and counit  are given by evaluation and coevaluation,\\
$\xy
{\ar^<<<<<<<{coev_{B}}(0,0)*+{B};(700,0)*+{\uno}};
\endxy$
and $\xy
{\ar^<<<<<<{ev_{B}}(0,0)*+{\iho \times \textnormal{A}};(700,0)*+{B}};
\endxy$.\\
That these lift directly to the case of \pb is the content of the following lemma.
\begin{lift}
\begin{enumerate}
\item The restriction of $\xy
{\ar^<<<<<<{ev_{B}}(0,0)*+{\iho \times \textnormal{A}};(700,0)*+{B}};
\endxy$ to $\ih{B} \times \textnormal{A}$ preserves pullbacks (lies in \pb).  We define the counit components via the restriction as  $\xy
{\ar^<<<<<<{ev_{B}}(0,0)*+{\ih{B} \times \textnormal{A}};(700,0)*+{B}};
\endxy$ and these are natural in B.
\item $\xy
{\ar^<<<<<<<{coev_{B}}(0,0)*+{B};(700,0)*+{\uno}};
\endxy$
preserves pullbacks, and its image lies in \un{B}.  We define the unit components via this factorization as $\xy
{\ar^<<<<<<<{coev_{B}}(0,0)*+{B};(700,0)*+{\un{B}}};
\endxy$ and these are natural in B.
\end{enumerate}

\end{lift}
\begin{proof}
\begin{enumerate}
\item
Let
$$\xy 
(0,0)*+{(P,a)}="p";(500,0)*+{(F,b)}="f"; (0,-400)*+{(G,c)}="g"; (500,-400)*+{(H,d)}="h";
{\ar@{=>}^{(r,\alpha)} "p";"f"}; {\ar@{=>}_{(s,\beta)} "p";"g"}; {\ar@{=>}^{(t,\theta)} "f";"h"}; {\ar@{=>}^{(u,\phi)}"g";"h"}; 
(50,-120)*{};(100,-120)*{} **\dir{-};
(100,-70)*{};(100,-120)*{} **\dir{-};
\endxy$$
be a pullback diagram in $\ih{B} \times \textnormal{A}$ (corresponding to a pullback square in each of \ih{B} and A).\\
The image of this pullback square under $ev_{B}$ is the outer square of
$$
\xy
(0,0)*+{Pa}="pa";(500,0)*+{Fa}="fa"; (1000,0)*+{Fb}="fb";
(0,-400)*+{Ga}="ga";(500,-400)*+{Ha}="ha"; (1000,-400)*+{Hb}="hb";
(0,-800)*+{Gc}="gc";(500,-800)*+{Hc}="hc"; (1000,-800)*+{Hd}="hd";
{\ar^{r_{a}} "pa";"fa"}; {\ar^{F\alpha} "fa";"fb"}; 
{\ar^{u_{a}} "ga";"ha"}; {\ar^{H\alpha} "ha";"hb"}; 
{\ar_{u_{c}} "gc";"hc"}; {\ar_{H\phi} "hc";"hd"}; 
{\ar_{s_{a}} "pa";"ga"}; {\ar_{G\beta} "ga";"gc"}; 
{\ar_{t_{a}} "fa";"ha"}; {\ar_{H\beta} "ha";"hc"}; 
{\ar^{t_{b}} "fb";"hb"}; {\ar^{H\theta} "hb";"hd"}; 
(50,-120)*{};(100,-120)*{} **\dir{-};
(100,-70)*{};(100,-120)*{} **\dir{-};
(550,-120)*{};(600,-120)*{} **\dir{-};
(600,-70)*{};(600,-120)*{} **\dir{-};
(50,-520)*{};(100,-520)*{} **\dir{-};
(100,-470)*{};(100,-520)*{} **\dir{-};
(550,-520)*{};(600,-520)*{} **\dir{-};
(600,-470)*{};(600,-520)*{} **\dir{-};
(1200,-800)*{.};
\endxy
$$
The top left square is a pullback as it is a component of the pullback square in \ih{B}.  The bottom right square is a pullback as it is the image of the pullback square in A, under the pullback preserving functor H.  The top right and bottom left squares are pullbacks as both $t$ and $u$ are cartesian.  Consequently the outer square is a pullback square as desired.  It is straightforward to verify that the counit components so defined constitute a natural transformation $\xy
{\ar^<<<<<<{ev}(0,0)*+{\ih{-} \times \textnormal{A}};(700,0)*+{1_{\pb}}};
\endxy$.
\item To see that $\xy
{\ar^<<<<<<<<{coev_{B}}(0,0)*+{B};(700,0)*+{\uno}};
\endxy$
preserves pullbacks, note that we have the product in Cat $\uno \cong \textnormal{[A,B]} \times \textnormal{[A,A]}$, and that the following diagram commutes:
$$
\xy
(0,0)*+{\textnormal{B}}="a";(800,0)*+{\uno}="bab";(1100,300)*+{\textnormal{[A,B]}}="ba";(1100,-300)*+{\textnormal{[A,A]}}="bb";
{\ar|<<<<<<<<<<<{coev_{B}} "a";"bab"};{\ar^{} "bab";"ba"};
{\ar^{} "bab";"bb"};
{\ar_{\widehat{1_{A}}} "a";"bb"};
{\ar^{\Delta} "a";"ba"};
\endxy
$$ where $\Delta(b)$ is the constant functor at $b$ for each object $b$ of B, $\widehat{1_{A}}(b) = 1_{A}$ is the identity functor on A, and the unlabelled arrows are the projections from the product.  Both $\Delta$ and $\widehat{1_{A}}$ clearly preserve pullbacks, so $coev_{B} = (\Delta,\widehat{1_{A}})$ preserves pullbacks.\\
To see that the image of $coev_{B}$ lies in \un{B}, we must show firstly that given an object $b$ of B, the functor $coev_{B}(b) = (\Delta,\widehat{1_{A}})(b) = (\Delta(b),1_{A})$ preserves pullbacks.  Certainly the constant functor $\Delta(b)$ at $b$ preserves pullbacks, as does $1_{A}$, so that $coev_{B}(b) = (\Delta(b),1_{A})$ preserves pullbacks.
Given a morphism $\xy
{\ar^{\alpha}(0,0)*+{a};(300,0)*+{b}};
\endxy$of B, we must verify that the natural transformation $coev_{B}(\alpha)$ is cartesian.  Now $coev_{B}(\alpha) = (\Delta(\alpha),1_{1_{A}})$, and as both $\Delta(\alpha)$ and $1_{1_{A}}$ are cartesian, it follows that  $coev_{B}(\alpha)$ is.\\
It is straightforward to verify that the unit components so defined constitute a natural transformation
$\xy
{\ar^<<<<<{coev}(0,0)*+{1_{\pb}};(700,0)*+{\un{-}}};
\endxy$.
\end{enumerate}
\end{proof}
\begin{closed}
The internal hom, unit and counit defined thus far give \pb the structure of a cartesian closed category.
\end{closed}
\begin{proof}
It remains to verify the triangle equations for the unit and counit.  Being defined exactly as in the case of Cat (where the triangle equations hold) they certainly hold in \pb.  Therefore \pb is cartesian closed.
\end{proof}
As \pb is cartesian closed we may enrich over it.   \pb obtains the structure of a \pb-category, \bpb, by defining $\bpb \textnormal{(A,B)} =\ih{B}$ for A and B categories with pullbacks.  Indeed \bpb is a cartesian closed \pb-category.\\
Now every \pb-category in particular has an underlying Cat-category (2-category).  To be precise (following \cite{x}) the forgetful functor $$\xy
{\ar^{\textnormal{U}}(0,0)*+{\pb};(500,0)*+{\textnormal{Cat}}};
\endxy$$ is finite product preserving and thus induces a 2-functor
$$\xy
{\ar^{\textnormal{U*}}(0,0)*+{\vcat{\pb}};(700,0)*+{\vcat{Cat}}};
\endxy \hspace{0.1cm}.$$
Given a \pb-category A, U*A has objects as A, and for objects $a$,$b$ of A, U*A($a$,$b$) = U(A($a$,$b$)).\\  In particular U*\bpb is the 2-category consisting of categories with pullbacks, pullback preserving functors and cartesian natural transformations.  We have seen that the underlying category, \pb, of U*\bpb is cartesian closed.  We conclude by showing that U*\bpb is a cartesian closed 2-category.
\begin{closed2cat}
The 2-category of categories with pullbacks, pullback preserving functors and cartesian natural transformations, U*\bpb, is a cartesian closed 2-category.
\end{closed2cat}
\begin{proof}
For A,B,C objects of \pb we have $$\textnormal{U*}\bpb(\textnormal{A} \times \textnormal{B,C)} =\textnormal{U}(\bpb(\textnormal{A} \times \textnormal{B,C)}) = \textnormal{U}(\textnormal{[A} \times \textnormal{B,C]}_{\textnormal{pb}}) \cong \textnormal{U}(\textnormal{[A,[B,C]}_{\textnormal{pb}}\textnormal{]}_{\textnormal{pb}}) = \textnormal{U*\bpb}(\textnormal{A,[B,C}]_{\textnormal{pb}})$$ naturally in A and C.
\end{proof}
I would like to thank my supervisor Dr. Stephen Lack for asking me whether \pb is cartesian closed.
\bibliographystyle{plain}

\end{document}